\documentclass[leqno]{amsart}

\usepackage[ansinew]{inputenc}
\usepackage{ngerman}
\usepackage{amsmath}
\usepackage{amssymb}
\usepackage{latexsym}
\usepackage{amsthm}
\usepackage{array}
\usepackage{booktabs}

\usepackage[usenames]{color}
\usepackage{endnotes} 
\let\footnote=\endnote

\DeclareMathOperator{\spa}{span}

\newtheorem{satz}{Satz}
\newtheorem*{satz-mercer}{Satz von Mercer}

\newtheorem*{Vermutete-Spurformel}{Vermutete Spurformel}

\newcommand{\R}{\mathbb{R}}

\newcommand{\N}{\mathbb{N}}

\begin{document}
 \title{Über Eigenwerte, Integrale und $\frac{\pi^2}{6}$: Die Idee der Spurformel}
\author{Daniel Grieser}\thanks{
Institut für Mathematik, Carl von Ossietzky Universität Oldenburg, D-26111 Oldenburg, e-mail: grieser@mathematik.uni-oldenburg.de}

\maketitle
\begin{quote}
\begin{center}\bf Zusammenfassung \end{center}

 Ausgehend von der Tatsache, dass die Summe der Diagonal-Elemente einer quadratischen Matrix gleich der Summe ihrer Eigenwerte ist, wird durch Analogie-Bildung eine analoge Formel ('Spurformel') für stetige symmetrische Funktionen \(G:[0,1]\times [0,1]\to\R\) geraten. Aus dieser wird dann die Identität \(\sum_{k=1}^\infty \frac1{k^2} = \frac{\pi^2}6\) hergeleitet. Konsequente  Fortführung der Analogie führt zu einer Beweisskizze der Spurformel. Es folgen weitere interessante Anwendungen der Spurformel.
 Dies ist die Ausarbeitung eines Vortrags, der sich an Mathematik-Studierende im zweiten Semester richtete und neben der Heranführung an spannende Mathematik zum Ziel hatte, zu zeigen, dass formale Analogiebildung und Raten fruchtbare Prozesse in der Mathematik sein können.
\end{quote}
\section*{Vorbemerkung}
Spurformeln sind eines der wichtigen Werkzeuge der modernen Mathematik. Sie verbinden so weit auseinanderliegende Bereiche wie Zahlentheorie, dynamische Systeme und partielle Differentialgleichungen. Sie helfen bei der Untersuchung der schwierigen Probleme der inversen Spektraltheorie ('Kann man die Form einer Trommel hören?'), und nahe verwandte Ideen führen zu einem Beweis des berühmten Indexsatzes von Atiyah und Singer und zu weitreichenden Verallgemeinerungen und Verfeinerungen
dieses Satzes. In der Physik findet die Welle-Teilchen-Dualität und die Korrespondenz von klassischer und Quantenmechanik Ausdruck in einer Spurformel.\footnote{Inverse Spektraltheorie und die Beziehung zur Physik werden in diesem Artikel erklärt. Für den Indexsatz von Atiyah und Singer verweisen wir die Leserin auf das gut lesbare Buch \cite{BooBle:TAASIFGTP}.}
Gleichzeitig ist die Grundidee der Spurformel schon mit sehr geringem Vorwissen, das in den ersten zwei Semestern eines Mathematikstudiums erworben wird, zu verstehen, und führt zu so erstaunlichen Ergebnissen wie
\begin{equation}
\label{pi6}
\sum\limits_{k=1}^\infty \frac1{k^2} = \frac{\pi^2}6.
\end{equation}
Dieser Artikel ist die  Ausarbeitung eines Vortrags, den ich für interessierte Studierende des zweiten Semesters im Sommersemester 2006 an der Carl von Ossietzky Universität Oldenburg gehalten habe. Ausgehend von einer einfachen bekannten Aussage der linearen Algebra 'leite' ich die Spurformel durch eine formale Analogiebildung 'her'. Die Analogiebildung besteht  im Wesentlichen im Ersetzen von Summen durch Integrale. Bevor ich nachprüfe, ob das Ergebnis stimmt, zeige ich, dass bei diesem Prozess wirklich etwas Neues, Kraftvolles, entstanden ist, indem ich die Formel \eqref{pi6} herleite. Es folgt ein nicht ganz vollständiger Beweis der Spurformel. Auch dieser ist, bei konsequenter Fortführung der Analogie, sehr einfach. Die fehlenden Teile sind gewisse Tatsachen aus der Funktionalanalysis (oder, im betrachteten Spezialfall, über Fourierreihen), deren vollständiger Beweis zwar aufwändiger ist, die aber
durch ihre Ähnlichkeit zu Aussagen der linearen Algebra nicht überraschen.
Schließlich erkläre ich, wie die Anwendung der Spurformel in  anderen Kontexten zu sehr tiefen Beziehungen wie den eingangs erwähnten führt.

Die Bildung gewagter Analogien, das Spielen mit und Modifizieren von Konzepten, ist ein höchst fruchtbarer Prozess in der Mathematik. Diese Tatsache sollte Studierenden der Mathematik, die meist nur den in der Lehre üblichen logischen deduktiven Aufbau kennenlernen, frühzeitig nahegebracht werden.\footnote{Nicht nur Studierenden, sondern auch Schülern! Die Spurformel ist hierfür natürlich weniger geeignet, aber es gibt auch für Schüler unzählige Themen, an denen dies durchgeführt werden kann. Siehe hierzu auch das wunderbare, schon für Schüler lesbare Buch \cite{Pol:MPRIAM} von Polya.}
Die Spurformel bildet hierfür ein hervorragendes Beispiel, gleichzeitig zeigt sie die Einheit der Mathematik und öffnet die Tür zur (sehr viel) höheren Mathematik.

Die Zielgruppe des Vortrags waren Mathematik-Studierende im 2. Semester, die Analysis I, Lineare Algebra I und in Analysis II ein wenig über Differentialgleichungen gelernt hatten. Ein zusätzlicher Reiz des Themas besteht darin, dass es auch für gestandene Mathematiker und Physiker interessant ist. Am Ende dieser Ausarbeitung und in Fußnoten, die ebenfalls am Ende gesammelt sind, sind daher ein paar weiterführende Bemerkungen und Literaturhinweise angefügt.
Insgesamt behalte ich den etwas informellen Stil des Vortrags weitgehend bei.

\section{Einleitung}
Dass Summen und Integrale eng verwandt sind, wird niemanden überraschen,
der sich einmal über die Bedeutung von Flächeninhalten unter Funktionsgraphen
und ihre Approximation durch Rechtecke (``Riemannsche Summen'')
Gedanken gemacht hat.

Daher sollte man ruhig einmal versuchen, in Formeln, die Summationen
enthalten, diese durch Integrale zu ersetzen (oder umgekehrt) -- vielleicht
kommt ja etwas Interessantes heraus. Spielerisches Erkunden führt oft zu
Erkenntnis! In der Analysis I haben wir schon einige Beispiele hierfür kennengelernt, z.B. des Integralkriterium für die Konvergenz von Reihen (z.B. \(\sum_{k=1}^\infty \frac1k\) divergiert, weil \(\int_1^\infty \frac1x\,dx=\lim_{N\to\infty} \log N\) divergiert).

Meist wird man erwarten, dass bei diesem einfachen Ersetzen etwas
schief geht -- aus einer exakten Formel wird eine nur ungefähr richtige. Zum Beispiel ist \(\int_1^\infty \frac1 {x^2}\,dx = 1\), aber offenbar \(\sum_{k=1}^\infty \frac1 {k^2} > 1\).

Hier soll eine Formel über die Eigenwerte von Matrizen vorgestellt werden,
die bei einer derartigen Ersetzung wieder eine exakte Formel liefert. Und
diese Formel hat erstaunliche Anwendungen.

\section{Die Spurformel für Matrizen}

Erinnern wir uns an etwas lineare Algebra: Sei $A$ eine reelle
$n\times n$-Matrix. Ein \textbf{Eigenwert} von $A$ ist eine Zahl
$\lambda\in \R$, für die es einen Vektor $u\in \R^n$, $u\not= 0$ gibt mit
\begin{equation} \label{eq0}
   Au = \lambda u
\end{equation}
Man berechnet die Eigenwerte wie folgt: Sei $p(z) = \det(zI-A)$
das charakteristische Polynom von $A$.
Dann sind die Eigenwerte die Nullstellen von $p$.

 Zur Berechnung der
Eigenwerte muss man die Nullstellen eines Polynoms bestimmen, das ist
schwierig. Ihre Summe ist aber leicht direkt aus $A$ abzulesen:

\begin{satz}[Spurformel für Matrizen]
Seien $\lambda_1,\lambda_2,\ldots,\lambda_n$ die Eigenwerte von $A=
(A_{ij})_{i,j=1,\ldots,n}$, aufgezählt mit Multiplizität. Dann ist
\begin{equation} \label{eq1}
   \sum^n_{k=1} \lambda_k = \sum^n_{i=1} A_{ii}
\end{equation}
``Aufgezählt mit Multiplizität'' bedeutet:
Jedes $\lambda$ taucht so oft auf, wie es seine Vielfachheit als Nullstelle des charakteristischen Polynoms von \(A\) angibt.

Falls \(A\) symmetrisch ist, d.h.\ \(A_{ij} = A_{ji}\ \forall i,j\), so ist die Vielfachheit von \(\lambda\) gleich der Dimension des Eigenraums $\{u:
Au = \lambda u\}$.
\end{satz}
Man nennt $\sum\limits^n_{i=1} A_{ii}$ die
\textbf{Spur von $\pmb{A}$}, daher der Name Spurformel.
\begin{proof}
Dies lässt sich zum Beispiel so beweisen, dass man den Koeffizienten von \(z^{n-1}\) des charakteristischen Polynoms \(p(z)\) auf zwei Weisen berechnet.\footnote{Skizze: Einmal mittels der Linearfaktorzerlegung \(
   p(z) = (z-\lambda_1)(z-\lambda_2) \cdots (z-\lambda_n)\, ,
\)
das ergibt $-\sum\limits^n_{k=1} \lambda_k$, und einmal mittels der \glq Leibniz-Formel\grq\ für die Determinante (Summe über alle Permutationen $\ldots$):
Da ein Term mit $z^{n-1}$ nur von dem Summanden $(z-A_{11})
\cdot\ldots\cdot (z-A_{nn})$ stammen kann, ist der Koeffizient
von $z^{n-1}$ gleich
\(
   -\sum^n_{i=1} A_{ii}\, .
\), und damit folgt \eqref{eq1}.
Die Gleichheit von algebraischer und geometrischer Multiplizität für symmetrische Matrizen ist ein zentrales Ergebnis der linearen Algebra.
}

Man beachte jedoch, dass für symmetrische Matrizen die Aussage des Satzes gar nichts mit Determinanten zu tun hat! Zur Vorbereitung auf die spätere Verallgemeinerung, wo uns kein Determinanten-Begriff zur Verfügung steht, geben wir einen \glq determinantenfreien\grq\ Beweis, der aber so nur für symmetrische Matrizen funktioniert.

Wenn $A$ symmetrisch ist, gibt es eine Orthonormalbasis aus Eigenvektoren $u^1,\dots,u^n$ zu den Eigenwerten $\lambda_1,\dots,\lambda_n$. Da $u^1,\dots,u^n$ eine Basis ist, lässt sich jeder Vektor $v\in\R^n$  als
\begin{equation}
\label{linkomb}
v=\sum_{k=1}^n \alpha_k u^k
\end{equation}
schreiben, und die Koeffizienten $\alpha_k$ sind einfach zu bestimmen: Nimmt man das Skalarprodukt von \eqref{linkomb} mit $u^{k_0}$, so erhält man $(v,u^{k_0}) =
\sum_{k=1}^n \alpha_k (u^k,u^{k_0}) = \alpha_{k_0}$ wegen der Orthonormalität, d.h.
\begin{equation}
\label{alphak}
\alpha_k = (v,u^k)\qquad \text{ für alle }k.
\end{equation}
Sei nun $A_i$ der $i$-te Zeilenvektor von $A$. Wir wenden \eqref{linkomb}, \eqref{alphak} auf $v=A_i$ an und erhalten
\begin{equation}
\label{Ai}
A_i = \sum_{k=1}^n (A_i,u^k)\,u^k\qquad\text{ für alle } i.
\end{equation}
Nun ist offenbar $(A_i,u^k)$ gerade die $i$-te Komponente des Vektors $Au^k$, und mit $Au^k=\lambda_ku^k$ folgt $(A_i,u^k)=\lambda_k u^k_i$. Setzt man dies in \eqref{Ai} ein und nimmt in dieser Vektorgleichung auf beiden Seiten die $j$-te Komponente, erhält man
\begin{equation}
\label{Aij}
A_{ij} = \sum_{k=1}^n \lambda_k \,u_i^k\, u_j^k \qquad \text{ für alle }i,j.
\end{equation}
Setzt man nun $i=j$ und summiert über $i=1,\dots,n$ so erhält man wegen der Normalität $\sum_{i=1}^n (u_i^k)^2 =  |u^k|^2 = 1$ die Behauptung.
\end{proof}

\section{Von Summen zu Integralen: Die Spurformel}

Was passiert eigentlich, wenn man in \eqref{eq1} Summen durch Integrale
ersetzt? Probieren wir das, und sehen wir, ob man daraus etwas Sinnvolles
fabrizieren kann!

Zunächst die rechte Seite von \eqref{eq1}: Statt über $i\in \{1,\ldots,n\}$ zu summieren,
wollen wir integrieren, sagen wir über $x\in [0,1]$, der Einfachheit halber.
Statt mit Zahlen $A_{ij}$ für $i,j = 1,\ldots,n$ sollten wir also mit
Zahlen $A_{xy}$ für $x,y\in [0,1]$ beginnen, und die rechte Seite von \eqref{eq1}
könnten wir durch
\[
 \int^1_0 A_{xx} \, dx
\]
ersetzen. Reelle Zahlen $A_{xy}$ für alle $x,y\in [0,1]$ vorzugeben,
bedeutet nichts anderes, als eine Funktion $A: [0,1]\times [0,1]
\longrightarrow \R$ zu betrachten.

Was wird aus der linken Seite von \eqref{eq1}? Was soll ``Eigenwert'' für
so eine Funktion $A$ bedeuten?

Dazu müssen wir die Eigenwertgleichung \eqref{eq0} im neuen Kontext
interpretieren. Was ist $u$, was bedeutet $Au$?

$u\in \R^n$ is durch seine Komponenten $u_i$, $i = 1,\ldots,n$, gegeben.
Analog zu $A$ sollten wir dies durch $u_x$, $x\in [0,1]$, ersetzen. $Au$
ist der Vektor mit Komponenten $(Au)_i = \sum^n_{j=1} A_{ij} u_j$,
und auch dies hat eine unmittelbare Übersetzung:
\[
   (Au)_x = \int^1_0 A_{xy}\, u_y \, dy\, .
\]
Also: $u$ und $Au$  sind Funktionen $[0,1]\longrightarrow \R$, und $Au$ ist durch diese Formel gegeben.

Damit kann man wieder das Eigenwertproblem \eqref{eq0} formulieren: Für
welche $\lambda\in \R$ existiert eine Funktion $u: [0,1]\longrightarrow
\R$,
$u\not= 0$, mit $Au = \lambda u$?
Wir werden unten sehen, dass (unter gewissen Bedingungen an $A$) eine
\emph{unendliche Folge} von Eigenwerten $\lambda_1,\lambda_2,\ldots$
existiert.

Die Symmetrie von \(A\), also die Bedingung \(A_{ij} = A_{ji} \ \forall i,j\), sollte sicherlich durch die Bedingung \(A_{xy} = A_{yx} \ \forall x,y\) ersetzt werden.

In Tabelle~1 sind die eben gefundenen Analogien zusammengestellt. Statt $u$ und $A$ verwenden wir die üblicheren Buchstaben $f$ und $G$, und statt $f_x, G_{xy}$ die üblichere Notation $f(x), G(x,y)$.
(der Buchstabe $G$
kommt von der verbreiteten Bezeichnung ``Greensche Funktion''). Die Zeile 'Skalarprodukt' wird später erklärt.

Beachte: Der Buchstabe \(G\) wird sowohl für den Operator (d.h.\ die lineare Abbildung \(f\mapsto Gf\)) als auch für die Funktion \(G(x,y)\), den sogenannten Integralkern dieses Operators, verwendet. Dies entspricht der in der linearen Algebra üblichen Identifizierung einer linearen Abbildung mit der sie definierenden Matrix.

Damit können wir eine Vermutung
formulieren.
\medskip

\begin{Vermutete-Spurformel}
Gegeben sei eine Funktion \(G:[0,1]\times [0,1]\to\R\), die symmetrisch ist, d.h.
\(G(x,y)=G(y,x) \ \forall x,y\). \\
\(\lambda_1,\lambda_2,\dots\) seien diejenigen Zahlen, für die die Gleichung \(Gu=\lambda u\) Lösungen \(u:[0,1]\to\R\), \(u\not\equiv 0\) hat; jedes \(\lambda\) trete dabei in der Folge so  oft auf, wie die Dimension des Vektorraums \(\{u: Gu=\lambda u\}\) angibt.
Dann gilt
\begin{equation} \label{spurformel}
   \sum^\infty_{k=1} \lambda_k = \int^1_0 G(x,x) \, dx.
\end{equation}
\end{Vermutete-Spurformel}
Damit das Integral überhaupt existiert, brauchen wir natürlich gewisse  Annahmen an \(G\).

Wir werden sehen, dass z.B. für stetiges \(G\) mit positiven Eigenwerten die Formel wirklich stimmt, dass sie, angewendet auf spezielle \(G\), erstaunliche Ergebnisse liefert, und dass eine etwas allgemeinere Version (wobei das Intervall \([0,1]\) durch höherdimensionale Gebiete ersetzt wird)  hoch interessante Ergebnisse über Obertöne, Billiardkugelbahnen und mehr liefert.\footnote{Eine \glq allgemeinste\grq\ Version der Spurformel gibt es nicht. Neben der eben erwähnten ($G$ stetig und symmetrisch, alle $\lambda_k\geq 0$) verdient die Spurformel von Lidskii Erwähnung ($G$ \glq Spurklasse\grq, nicht notwendig symmetrisch, siehe \cite{Sim:TITA}). In den Anwendungen im Abschnitt \ref{secobertoene} wird eine distributionenwertige Version verwendet -- siehe die Anmerkungen dort --, die aus der Lidski-Spurformel hergeleitet werden kann (siehe die dort zitierte Literatur). Die Grundidee ist aber immer dieselbe: Gleichung \eqref{spurformel}.}

\begin{table}
\begin{center}

\newcommand{\blue}{\textcolor{Black}}
\newcommand{\green}{\textcolor{Black}}
\newcommand{\red}{\textcolor{Black}}
\newcommand{\mahogany}{\textcolor{Black}}

\setlength{\extrarowheight}{10pt}
\thicklines

\begin{tabular}{c|c}
\blue{\textit{\textbf{ Vektoren und Summen}}} &
\green{\textit{\textbf{ Funktionen und Integrale}}} \\ \toprule
\blue{$i,j \in\{ 1,\ldots,n\}$} & \green{$x,y \in [0,1]$}\\
\blue{$u\in \R^n$} & \green{$f:[0,1] \longrightarrow \R$}\\
\blue{Komponente \quad $u_i$} & \green{Funktionswert \quad $f(x)$}\\
\toprule
\blue{$n\times n$-Matrix $A = (A_{ij})_{i,j = 1,\ldots,n}$} &
\green{Funktion $G: [0,1] \times [0,1] \longrightarrow \R$}\\
\blue{Matrixelement $A_{ij}$} & \green{Funktionswert $G(x,y)$}\\ \toprule
\blue{$Au \in \R^n$} & \green{$Gf: [0,1] \longrightarrow \R$}\\
\blue{ist der Vektor mit Komponenten} & \green{ist die Funktion mit
Werten}\\
\blue{$(Au)_i = \displaystyle\sum^n_{j=1} A_{ij} u_j$} &
\green{$(Gf)(x) = \displaystyle\int^1_0 G(x,y) \, f(y) \, dy$}\\[4ex]
\toprule
\multicolumn{2}{c}{\textbf{Eigenwertgleichung}}\\
\blue{$Au = \lambda u$} & \green{$Gf = \lambda f$}\\[0ex]
\toprule
\multicolumn{2}{c}{\mahogany{\textbf{Spurformel}}}\\
\red{$\displaystyle\sum^n_{k=1} \lambda_k = \sum^n_{i=1} A_{ii}$} &
\red{$\displaystyle \sum^\infty_{k=1}\lambda_k = \int^1_0 G(x,x)\, dx$}
\\[8ex]
\toprule
\multicolumn{2}{c}{\textbf{Symmetrie}}\\
\blue{$A_{ij}=A_{ji}\quad \forall i,j$} &
\green{$G(x,y) = G(y,x)\quad \forall x,y$}
\\[4ex]
\multicolumn{2}{c}{\textbf{Skalarprodukt}}\\
\blue{$(u,v) = \displaystyle\sum_{i=1}^n u_i v_i$} &
\green {$(f,g) = \displaystyle\int_0^1 f(x) g(x)\, dx$}
\end{tabular}
\end{center}

\caption{Lexikon}
\end{table}

\section{$\frac{\pi^2}6$, Teil I: Ein einfaches Randwertproblem}

Bevor wir die vermutete Spurformel \eqref{spurformel}
beweisen, wollen wir sehen, dass sie interessant ist, indem wir sie auf ein
spezielles $G$ anwenden.

Die Funktion $G$, die wir verwenden wollen, tritt im Kontext einer
einfachen Differentialgleichung auf. Da diese Gleichung einfacher hinzuschreiben ist als
die Formel für $G$ und da wir sie später noch einmal brauchen, fangen wir
damit an.

Wir betrachten eine der einfachsten Differentialgleichungen, die es gibt (ein Randwertproblem):
\begin{align}
   -u'' (x) = f(x) \quad  \text{ für } x\in [0,1],\qquad   u(0)  =u(1) = 0\, .  \label{eq4}
\end{align}
Die Funktion $f$ ist gegeben, $u$ ist gesucht. Wie findet man $u$? Ganz einfach:
Die Gleichung \(-u''=f\) zweimal integrieren, dabei erhält man zwei freie
Integrationskonstanten. Diese wählt man so, dass die Randbedingung \(u(0)=u(1)=0\) erfüllt ist.

Hier ist das Ergebnis in einer Form, wie wir sie hier brauchen:

\begin{satz}\label{satz2}
Das Problem \eqref{eq4} hat zu jedem stetigen $f: [0,1]
\longrightarrow \R$ eine eindeutige Lösung $u: [0,1] \longrightarrow \R$.
Diese ist gegeben durch
\begin{align}
   u(x) & = \int^1_0 G(x,y) \, f(y) \, dy \label{eqloesung}
\intertext{mit}
   G(x,y) & = \begin{cases} x\cdot (1-y) &  \text{ für }x \leq y\\
              y \cdot (1-x)  & \text{ für } x\geq y \end{cases} \label{eq5}
\end{align}
\end{satz}

\begin{proof}[Beweis (alternativ zur skizzierten Herleitung)] \small
Man rechnet leicht nach, dass dieses $u$, anders geschrieben als
\[
   u(x) = (1-x) \int^x_0 y\, f(y)\, dy + x \int^1_x  (1-y)\, f(y)\, dy
\]
wirklich eine Lösung von \eqref{eq4} ist. Es bleibt die
Eindeutigkeit zu zeigen. Sind $u_1, u_2$ zwei Lösungen von \eqref{eq10},
so gilt für die Differenz $u= u_1-u_2$:
\begin{align*}
   u''(x)  = 0 \quad \forall x \in [0,1], \qquad
   u(0) = u(1)  = 0 \, .
\end{align*}
Das heißt,  $u$ ist eine lineare Funktion, die an zwei Punkten
verschwindet. Also muss $u \equiv 0$ sein, also $u_1 \equiv u_2$.
\end{proof}
Beachte, dass $G$ symmetrisch ist.

\section{$\frac{\pi^2}6$, Teil II: Berechnung der Eigenwerte}
\label{secewberechnung}

Wir wollen nun die Spurformel \eqref{spurformel} auf die Funktion \(G\) in \eqref{eq5} anwenden.
Hierzu müssen wir die Eigenwerte
$\lambda_k$ von $G$ bestimmen, also die Zahlen $\lambda$, für die es eine
stetige Funktion $f \not\equiv 0$ gibt mit
\begin{equation} \label{eq6}
   Gf = \lambda f\, .
\end{equation}
Wie löst man diese Gleichung? Verwendet man direkt \eqref{eq5}, so ist unklar, wie man diese Integralgleichung lösen soll. Es stellt sich heraus, dass es einfacher ist, das analoge Problem für die Differentialgleichung zu lösen, aus der wir \(G\)
erhalten hatten. Hierfür kann man \eqref{eq6} zweimal ableiten. Oder einfacher:   Schreibt man $u= Gf$, so löst \(u\) nach Satz \ref{satz2} das Randwertproblem
\eqref{eq4}, also wird \eqref{eq6} zu
\begin{align}
   u  = -\lambda u'',\qquad   u(0)  =u(1) = 0\, , \label{eq8}
\end{align}
und aus \(f\not\equiv 0\) folgt \(u\not\equiv 0\) wegen \(f=-u''\). Erfüllt $u\not\equiv 0$ umgekehrt \eqref{eq8}, so
folgt \eqref{eq6} für $f= -u''$, und wegen der Eindeutigkeit
in Satz \ref{satz2} muss $f\not\equiv 0$ sein. Also müssen wir die Zahlen $\lambda$
bestimmen, für die \eqref{eq8} eine Lösung $u\not\equiv 0$
hat. Das ist einfach:

\begin{satz}
Das Problem \eqref{eq8} hat eine Lösung $u\not\equiv 0$ genau
dann, wenn $\lambda = \frac{1}{\pi^2 k^2}$ für ein $k\in \N$ ist. In diesem
Fall ist
\begin{equation} \label{eq9}
   u(x) =a\cdot \sin k\pi x \qquad\qquad (a\in \R)\, .
\end{equation}
\end{satz}

\begin{proof}
Wegen $u\not\equiv 0$ muss $\lambda\not= 0$ sein. Die allgemeine Lösung der
Gleichung $u'' + \frac{1}{\lambda}\; u = 0$ ist wohlbekannt, nämlich
\begin{align*}
   u(x)  & = a \sin \frac{x}{\sqrt{\lambda}} + b \cos
             \frac{x}{\sqrt{\lambda}} \qquad \text{ für } \lambda >
0\\[-3ex]
\intertext{und}
   u(x)  & = a \, e^{-{x}/{\sqrt{-\lambda}}} + b \,
             e^{{x}/{\sqrt{-\lambda}}} \qquad\quad \text{ für }
            \lambda <  0\, ,
\end{align*}
mit beliebigen $a,b \in \R$.

Man sieht leicht, dass die Bedingung $u(0) = u(1) =0$ nur den Fall $\lambda
> 0$ zulässt, dass dann $b=0$ sein muss (wegen $u(0) = 0$), $a\neq 0
$ beliebig sein kann und $\frac{1}{\sqrt{\lambda}}$ ein ganzzahliges
Vielfaches von $\pi$ sein muss (wegen $u(1) = 0$).
\end{proof}

Wir haben also $\lambda_k = \frac{1}{\pi^2 k^2}$, $k\in \N$.
Die Spurformel \eqref{spurformel} (wenn sie denn stimmt!) ergibt nun:
\begin{align*}
   \sum^\infty_{k=1} \frac{1}{\pi^2 k^2}\;
 = \int^1_0 x\cdot(1-x)\, dx
 = \frac{1}{2} - \frac{1}{3} = \frac{1}{6}\, ,
\end{align*}
also
\begin{center}
\setlength{\fboxrule}{1pt}\setlength{\fboxsep}{4mm}
\fbox{$\displaystyle\sum^\infty_{k=1} \;\frac{1}{k^2} =
\frac{\pi^2}{6}$}
\end{center}

Falls unsere Analogiebildung sich präzise begründen lässt, haben wir uns also ein Werkzeug geschaffen, mit dem wir höchst nicht-triviale Dinge beweisen können!

\section{Stimmt die Spurformel?}

Eine formale Analogie hatte uns von der Matrix-Spurformel \eqref{eq1} auf die vermutete Integral-Spurformel \eqref{spurformel} geführt. Wir werden nun sehen, dass sich auch der {\em Beweis} der Matrix-Spurformel mittels derselben Analogie fast wörtlich in einen (Fast-)Beweis der Spurformel \eqref{spurformel} übersetzen lässt.
Wir setzen voraus, dass $G$ stetig ist.

Als ersten Schritt im Beweis der Matrix-Spurformel hatten wir den Satz der linearen Algebra zitiert, dass für jede symmetrische Matrix eine Orthonormalbasis aus Eigenvektoren existiert. Gilt ähnliches auch für den Operator $G$?

Zunächst: Was bedeutet \glq Orthonormalbasis\grq? Als erstes brauchen wir einen Vektorraum. Nach den Überlegungen oben sollten wir
\begin{equation} 
   V:= \{\text{stetige Funktionen } [0,1] \longrightarrow \R\}
\end{equation}
nehmen.
Um von Orthogonalität sprechen zu können, braucht man ein
\emph{Skalarprodukt}. Das Skalarprodukt von $u,v\in \R^n$ ist $(u,v) =
\sum^n_{i=1} u_i\, v_i$. In Fortsetzung unserer Analogie sollte das
Skalarprodukt zweier Funktionen $f,g:[0,1] \longrightarrow \R$ also die
Zahl $$(f,g) := \int^1_0 f(x)\, g(x)\, dx$$ sein.

Dass damit die aus der linearen Algebra im $\R^n$ bekannten Dinge stimmen,
bestärkt uns darin, dass dies eine sinnvolle Definition ist: Die Abbildung
 $f,g\mapsto (f,g)$ ist wirklich ein Skalarprodukt auf \(V\) (also bilinear,
symmetrisch und positiv definit), insbesondere ist
\begin{equation} \label{eq11}
   \|f\| := \sqrt{(f,f)} = \sqrt{\int^1_0 f(x)^2\, dx}
\end{equation}
eine Norm auf \(V\).\footnote{
Außerdem ist, wie in der linearen Algebra, die Symmetrie von \(G\) (also \eqref{eq10})
äquivalent zu
\[
   (Gf,g) = (f,Gg) \qquad \forall f,g \in V.
\]
Denn ausgeschrieben ist dies \(\int_0^1\int_0^1 G(x,y)f(y)g(x)\,dx\,dy = \int_0^1\int_0^1 G(y,x) f(y)g(x)\,dx\,dy\) -- wobei hier die Integrationsreihenfolgen vertauscht wurden --, und dies kann nur dann für alle \(f,g\) stimmen, wenn \(G(x,y)=G(y,x)\ \forall x,y\) gilt. }

Eine Orthonormalbasis ist nun eine Folge von Funktionen $f_k\in V$ mit:
\begin{itemize}
\item
$(f_k,f_l) = 0 \quad \forall k\not= l$ (orthogonal)
\item
$\|f_k\| = 1 \quad \forall k$ (normiert)
\item
$\forall f \in V \quad \exists \alpha_1,\alpha_2,\ldots \in \R$ mit
\end{itemize}
\begin{equation} \label{eq12}
   f = \sum^{\infty}_{k=1} \alpha_k f_k \qquad \qquad \text{(Basis)}
\end{equation}

Vorsicht: Dieser Basisbegriff ist ähnlich, aber nicht gleich dem der linearen
Algebra.
Dort müssen alle Summen endlich sein! Bei unendlichen Summen muss man sich
über Konvergenz Gedanken machen. \eqref{eq12} soll also bedeuten:
\begin{equation} \label{eq13}
   \|f - \sum^N_{k=1} \alpha_k f_k \| \longrightarrow 0 \qquad (N\to
                                                              \infty)
\end{equation}
mit der Norm \eqref{eq11}, also ganz analog zu Reihen in $\R$.

Es gilt nun wirklich:
\begin{satz} \label{satzspektralsatz}
Sei $G: [0,1] \times [0,1]\longrightarrow \R$ stetig und symmetrisch, d.h.:
\begin{equation} \label{eq10}
   G(x,y) = G(y,x) \qquad \forall x,y \in [0,1]
\end{equation}
Dann gibt es eine Zahlenfolge $\lambda_1,\lambda_2,\ldots$  und eine Orthonormalbasis zugehöriger
Eigenfunktionen
$f_k: [0,1]\longrightarrow \R$:
\begin{equation}
   Gf_k = \lambda_k\, f_k \quad \forall k\label{fkeigenfkt}
\end{equation}
Alle $f_k$ sind stetig, und $\lambda_k \to 0$ für $k\to \infty$.
\end{satz}

Man beweist Satz \ref{satzspektralsatz} in einer einführenden
Funktional\-analy\-sis-Vorlesung, er ist ein Spezialfall des
``Spektralsatzes für kompakte Operatoren'', siehe z.B. \cite{Wer:F}, Theorem VI.3.2 und Kapitel VI.4.
In unserem Beispiel sind die $f_k(x) =
\sqrt{2}
\sin  k\pi x$, und der Satz drückt eine grundlegende Aussage über
Fourierreihen aus (sogenannte \(L^2\)-Konvergenz).

\begin{proof}[Beweis der Spurformel]
Wir imitieren nun den Beweis der Matrix-Spurformel \eqref{eq1}. Hierbei vernachlässigen wir zunächst Konvergenzfragen.

Indem man \eqref{eq12} skalar mit $f_{k_0}$ multipliziert und dann $k_0$ durch $k$ ersetzt, erhält man analog zu \eqref{alphak}
\begin{equation} \label{eq14}
   \alpha_k = (f,f_k) \qquad \text{ für alle } k\, .
\end{equation}
Wir wenden nun \eqref{eq12}, \eqref{eq14} auf die Funktion $f(y) = G(x,y)$
an, für ein festes $x\in [0,1]$, und erhalten analog zu \eqref{Ai}
\begin{equation} \label{eq15}
   G(x,y) = \sum_k (G(x,\cdot),f_k) \; f_k(y)
\end{equation}
Nun ist
\(   (G(x,\cdot),f_k)  = \int_0^1 G(x,y)\, f_k(y)\, dy  = (Gf_k)(x)\), und wegen $Gf_k = \lambda_k f_k$ ist dies gleich $\lambda_k f_k(x)$. Also wird \eqref{eq15} zu
\begin{equation} \label{eq16}
   G(x,y) = \sum^\infty_{k=1} \lambda_k\, f_k(x) \, f_k(y)\, ,
\end{equation}
analog zu \eqref{Aij}.
Wir setzen nun $x=y$ und integrieren:
\begin{align}
    \int^1_0 G(x,x)dx & = \int^1_0 \sum^\infty_{k=1} \lambda_k\, f_k(x)^2\,
                         dx  \label{eq17}\\[1ex]
                      & = \sum^\infty_{k=1} \lambda_k \int^1_0 f_k(x)^2\,
                         dx \label{eq18}\\[1ex]
                      & = \sum^\infty_{k=1} \lambda_k \, ,
\end{align}
da $1 = \|f_k\|^2 = \int^1_0 f_k(x)^2\, dx$. Das ist die Spurformel!
\end{proof}

Dafür, dass man damit so tolle Dinge wie $\sum_k \frac{1}{k^2} =
\frac{\pi^2}{6}$ beweisen kann, ist es ein erstaunlich einfacher Beweis.

War das ein Beweis? Wie steht es mit Konvergenz etc.?
\begin{itemize}
\item
Für \eqref{eq14} wurde die unendliche Summe aus dem Skalarprodukt
herausgezogen. Das ist leicht zu rechtfertigen (Stetigkeit des Skalarprodukts).
\item
Problematischer wird es bei \eqref{eq15}: Satz \ref{satzspektralsatz} garantiert nur, dass die
Konvergenz $\sum^N_{k=1} \alpha_k f_k \xrightarrow{N\to \infty} f$ in der Norm
$\|\;\|$, also ``im quadratischen Mittel'', gilt. Sie stimmt im Allgemeinen
nicht punktweise, d.~h. es kann einzelne $y \in [0,1]$ geben mit
\(
   \sum^\infty_{k=1} \alpha_k\, f_k(y) \not= f(y)
\), oder wo die Summe sogar divergiert.
Daher braucht \eqref{eq15} und damit \eqref{eq16} nicht punktweise (für jedes $y$) zu stimmen.
\item
In \eqref{eq18} haben wir Integral und unendliche Summe vertauscht. Dies ist nur unter gewissen Voraussetzungen, z.B.\ bei gleichmäßiger Konvergenz, erlaubt.
\end{itemize}
Der zweite und dritte Punkt lassen sich mit etwas mehr Arbeit rechtfertigen. Es gilt nämlich:
\begin{satz-mercer}
Falls alle $\lambda_k \geq 0$ sind, so gilt \eqref{eq16} für alle $x,y \in
[0,1]$, und die Konvergenz ist gleichmäßig.
\end{satz-mercer}
(Siehe zum Beispiel \cite{Wer:F}, Satz VI.4.2).
In unserem Beispiel ist \(f_k(x)= \sqrt 2 \sin k\pi x\) und \(\lambda_k = \frac 1 {\pi^2 k^2}\). In diesem Fall braucht man nicht auf Satz \ref{satzspektralsatz} und den Satz von Mercer zu verweisen, die Aussagen folgen aus Sätzen über Fourierreihen.%
\footnote{Wegen der gleichmäßigen Beschränktheit der \(f_k\) und der Konvergenz der Reihe \(\sum_k \lambda_k\) konvergiert die rechte Seite von \eqref{eq16} gleichmäßig. Der Satz von Fej\'er ergibt dann, dass die Summe wirklich \(G(x,y)\) ist, da diese Funktion stetig ist.
}

\section{Ein neuer Blickwinkel} \label{secneublick}
Obwohl die Spurformel eine Aussage über Integrale und Summen ist, war in unserem Beispiel das zentrale Objekt die Differentialgleichung \eqref{eq4}. Denn nur mit deren Hilfe gelang es uns, die Eigenwerte \(\lambda_k\) zu bestimmen.

Daher liegt es nahe, andere Differentialgleichungen zu betrachten, um weitere interessante Resultate zu erhalten. Das Grundprinzip ist dabei folgendes.

Wir betrachten eine lineare Differentialgleichung
\begin{equation}\label{DGL}
 P(u) = f, \quad \text{Randbedingung an }u,
\end{equation}
(in unserem Beispiel ist \(P(u)=-u''\) und die Randbedingung ist \(u(0)=u(1)=0\)).
Eine zugehörige Spurformel ergibt sich, wenn man zwei grundverschiedene Ansätze zur Lösung dieser Gleichung verfolgt und die Resultate gleichsetzt\footnote{Damit das funktioniert, müssen gewisse Annahmen an \(P\) gemacht werden, z.B. Selbstadjungiertheit, Invertierbarkeit und Diskretheit des Spektrums. In allen hier betrachteten Beispielen sind diese Bedingungen erfüllt. Die sich ergebende Spurformel ist dann, genau genommen, die für den Operator \(P^{-1}\).

Die im Text folgende Überlegung gibt eine alternative Herleitung der
fundamentalen Identität \eqref{eq16} und damit der Spurformel (im Fall \(G=P^{-1}\)). Eine weitere Herleitung ist folgende:
Ist $P_k$ die orthogonale Projektion auf den Unterraum $\spa\{f_k\}$, so
hat $P_k$ den Integralkern $P_k(x,y) = f_k(x)\, f_k(y)$, d.h. \((P_kf)(x) = \int P_k(x,y)f(y)\,dy\ \forall f\) (einfach
nachrechnen). Satz \ref{satzspektralsatz} lässt sich umformulieren als $G = \sum_k \lambda_k\,
P_k$. Nimmt man den Integralkern beider Seiten, folgt $G(x,y)=\sum_k \lambda_k\, f_k(x)\, f_k(y)$.

Dies lässt sich übrigens leicht in einen weiteren Beweis der Spurformel für diagonalisierbare Matrizen übersetzen.
}.

\begin{itemize}
\item[1.] \emph{Lösung mittels der Spektralzerlegung}\\[1ex]
Wir lösen das Eigenwertproblem für \eqref{DGL}, d.h. wir bestimmen Zahlen \(\mu_k\) und eine Orthonormalbasis \((f_k)\) mit  $Pf_k = \mu_k\, f_k$ für alle \(k\) (und so, dass alle \(f_k\) die Randbedingung erfüllen). Ist nun in \eqref{DGL}
die rechte Seite $f=f_k$ für ein $k$, so hat man die Lösung $u =
\frac{1}{\mu_k}\; f_k$.

Für allgemeines $f$ schreiben wir
\[
   f= \sum_k \alpha_k \, f_k, \qquad \alpha_k = (f,f_k)\, ,
\]
und wegen Linearität (und Stetigkeit $\ldots$) ist dann
\[
   u = \sum_k \frac{1}{\mu_k}\; \alpha_k \, f_k
\]
eine Lösung von \eqref{DGL}.

Schreiben wir \(\lambda_k=\frac 1 {\mu_k}\), so können wir dies wegen $\alpha_k = (f,f_k) = \int^1_0 f(y)\, f_k(y)\, dy$ als
\begin{align}
   u(x) & = \sum_k \lambda_k \, f_k(x) \int^1_0 f(y)\, f_k(y)\, dy
            \notag\\[1ex]
        & = \int^1_0 \left(\sum_k \lambda_k \, f_k(x)\, f_k(y)\right) f(y)\,
            dy \label{eq20}
\end{align}
schreiben.
\item[2.] \emph{Lösung durch eine alternative, \glq direkte\grq\ Methode}\\[1ex]
Bei unserem Beispiel ist dies zweimaliges Integrieren und Bestimmen der Konstanten, dies ergibt
\begin{equation} \label{eq21}
   u(x) = \int^1_0 G(x,y)\, f(y)\, dy
\end{equation}
mit $G$ wie in \eqref{eq5}.
\end{itemize}
Da \eqref{eq20}, \eqref{eq21} für alle stetigen $f$ gelten, folgt $
\sum_k \lambda_k\, f_k(x)\, f_k(y)=G(x,y)$, und indem man $x=y$ setzt und
integriert, folgt $\sum_k \lambda_k = \int G(x,x)\,dx$.

\section{Obertöne und Billiardkugelbahnen} \label{secobertoene}

Wendet man die Idee der Spurformel auf andere lineare Gleichungen an, die
neben der spektralen eine direkte Lösung besitzen, erhält man mitunter
erstaunliche Ergebnisse. Dies sind manchmal Formeln und manchmal qualitative Aussagen. Als zweites Beispiel erläutern wir eine Beziehung zwischen
Obertönen und Billiardkugelbahnen einer ebenen Membran.

Statt des Intervalls $[0,1]$ betrachten wir ein Gebiet $\Omega \subset
\R^2$ (beschränkt, mit stückweise glattem Rand $\partial \Omega$, z.~B. die
Kreisscheibe oder ein Rechteck).

\subsection{Obertöne: Die Eigenwerte des Laplace-Operators}
Eine natürliche Verallgemeinerung des Randwertproblems \eqref{eq4} auf
zwei Dimensionen ist das sogenannte Dirichlet-Problem
\begin{align} \label{eq22}
   -\Delta u  & = f \qquad \text{ in } \Omega\\
   u(x)       & = 0 \qquad \text{ für alle } x \in \partial\Omega
\label{eq23}
\end{align}
für eine gegebene Funktion $f: \Omega \longrightarrow \R$ und gesuchtes $u:
\Omega \longrightarrow \R$: Hierbei ist
\[
   \Delta u :=  \frac{\partial^2 u}{\partial x^2_1} + \frac{\partial^2
    u}{\partial x^2_2} \qquad\text{ der \emph{Laplace-Operator.}}
\]
Gemäß der Strategie in Abschnitt \ref{secneublick} betrachten wir
das zugehörige Eigenwertproblem, wo \eqref{eq22} durch

\[
   -\Delta u = \mu u \qquad \text{ in } \Omega
\]
mit $\mu\in \R$ ersetzt ist.

Ähnlich zu Satz \ref{satzspektralsatz} (wo \(\lambda=\frac1\mu\) ist) kann man zeigen, dass es Folgen
\[
\begin{array}{rccccc}
0 & < &  \mu_1 & \leq & \mu_2 & \leq \ldots \qquad\longrightarrow
\infty\\
  &      &  u_1 &      & u_2   &
\end{array}
\]
von Eigenwert-Eigenfunktionspaaren derart gibt, dass die $u_k$ eine
Orthonormalbasis bilden. Z.~B. sind für das Quadrat $[0,1] \times [0,1]$
die Eigenwerte gerade die Zahlen $\pi^2 (n^2+m^2)$ mit $n,m\in \N$,
aufsteigend angeordnet, mit Eigenfunktionen $2 \sin \pi n x_1 \, \sin \pi m x_2$ (nachrechnen!\footnote{Dass dies wirklich alle Eigenfunktionen sind, ist nicht so einfach zu zeigen. Dies lässt sich zum Beispiel mit ein wenig Theorie der Fourierreihen beweisen.}). Für die meisten Gebiete $\Omega$ kann man die
Eigenwerte aber nicht explizit berechnen.

Wir haben nun folgende Situation: Jedem Gebiet $\Omega \subset \R^2$ ist
eine Zahlenfolge $0< \mu_1 \leq \mu_2 \leq \ldots \to \infty$
zugeordnet.

Dies legt viele Fragen nahe, z.~B.: Welche Zahlenfolgen kommen vor, kann man
aus den Eigenwerten das Gebiet rekonstruieren (sogenanntes inverses Spektralproblem), kann man die $\mu_k$ zumindest
``ungefähr'' durch geometrische Größen von $\Omega$ beschreiben?

Man stelle sich nun $\Omega$ als Trommel vor, die am Rand $\partial\Omega$ fest eingespannt ist. Der Klang, der beim Schlagen der Trommel entsteht, lässt sich wie alle Klänge in einen Grundton und Obertöne zerlegen. Im Unterschied zu den Klängen, die von Saiten- oder Blasinstrumenten erzeugt werden, sind die Obertöne nicht harmonisch, d.h.\ ihre Frequenzen sind nicht das 2-, 3-, 4- \dots fache der Frequenz des Grundtons. Aus physikalischen Schwingungsgesetzen folgt, dass die Frequenz des Grundtons genau der Eigenwert $\mu_1$ und die Frequenzen der Obertöne genau $\mu_2,\mu_3.\dots$ sind, bis auf einen Normalisierungsfaktor, der z.~B. vom
Material abhängt (dies stimmt nur in einer idealisierten linearen
Welt). Siehe z.B. \cite{CouHil:MMP}. Die zweite Frage lässt sich also so formulieren:
\begin{quote}
Kann man die Form einer Trommel hören?
\end{quote}
Diese Frage ist bis heute offen (es gibt Gegenbeispiele mit Ecken, aber
bisher keine mit glattem Rand, siehe \cite{GorWebWol:OCHSD}.\footnote{Genau genommen hört man
nicht nur die Frequenzen der Obertöne, sondern auch deren relative
Lautstärken. Auch zu dieser verschärften Form der Frage gibt es
Gegenbeispiele, siehe \cite{BusConDoy:SPID}.}) Man kann aber zeigen, dass man viele geometrische
Eigenschaften von $\Omega$ hören, d.~h. aus den Eigenwerten bestimmen kann,
z.~B. den Flächeninhalt von \(\Omega\), die Länge des Randes und die Längen
geschlossener Billiardkugelbahnen in \(\Omega\). All dies erhält man mit Hilfe der Spurformel! Wir diskutieren hier nur den  letzten Punkt und erklären zunächst, was gemeint ist.

\subsection{Billiardkugelbahnen}
Wir stellen uns nun unser Gebiet $\Omega$ als Billiardtisch vor. Wir legen
eine (unendlich kleine) Billiardkugel auf den Tisch und stoßen sie. Unter
Idealbedingungen (keine Reibung, kein Effet) rollt sie gerade zum Rand,
wird dort nach dem Reflektionsgesetz (Einfallswinkel = Ausfallswinkel)
\glq reflektiert\grq, rollt dann wieder geradeaus bis zur nächsten Randberührung
usw., unendlich lange.\footnote{Was passiert, wenn die Kugel genau in eine Ecke rollt? Die für unseren Kontext sinnvollste Antwort ist nicht einfach, daher wollen wir im Folgenden nur Bahnen betrachten, die niemals eine Ecke treffen.}

Je nachdem, wo und in welche Richtung man die Kugel abgestoßen hat, kann
es passieren, dass die Kugel nach einer Weile genau zum Ausgangspunkt
zurückkehrt. Dies ist natürlich eine große Ausnahme. Kehrt sie sogar in
derselben Richtung zurück, wie sie abgestoßen wurde, wird sie danach immer
wieder dieselbe Bahn durchlaufen. Man spricht dann von einer
\emph{geschlossenen Billiardkugelbahn} in $\Omega$. Man nennt
\[
   L(\Omega) := \{\text{die Längen geschlossener Billiardkugelbahnen in
                }\Omega\} \subset \R
\]
das \emph{Längenspektrum} von $\Omega$. Mit Länge ist die Strecke gemeint,
die die Kugel bis zur Wiederkehr zurückgelegt hat.

Gibt es überhaupt solche Bahnen, ist also $L(\Omega) \not= \emptyset$? Man
kann zeigen, dass es immer welche gibt\footnote{%
Jedenfalls im Fall, dass der Rand von $\Omega$ glatt ist, dass also $\Omega$ keine Ecken hat. Existieren Ecken, so ist unbekannt, ob dies stimmt. Es ist sogar ein ungelöstes Problem, ob für jeden dreieckigen Billiardtisch eine geschlossene Billiardkugelbahn existiert! Eine Referenz für diese Themen ist \cite{Tab:GB}.}, für glattes konvexes $\Omega$ kann man
sogar auf dem Rand für jedes $n\in \N$ Punkte $p_1,\ldots,p_n$ finden, so
dass der Weg der Kugel genau $p_1 \to p_2 \to \ldots \to p_n \to p_1\to\dots$ ist.

\bigskip\noindent
\textbf{Übung:} Bestimme die geschlossenen Billiardkugelbahnen auf einem
kreisförmigen
Billiardtisch und auf einem rechteckigen Billiardtisch.

\subsection{Die Spurformel für die Wellengleichung}

Wir wollen nun die Beziehung

\bigskip\noindent
\parbox[t]{5.3cm}{Eigenwerte des Laplace-Opera\-tors auf
$\Omega$}
$\quad\leftrightsquigarrow\quad$
\parbox[t]{5.3cm}{Längen geschlossener Billiardkugelbahnen in $\Omega$}
\\[2ex]
erklären. Wir werden sehen, dass die Eigenwerte das Längenspektrum
eindeutig bestimmen.\footnote{Dies gilt jedenfalls für ``generische''
$\Omega$, d.~h. bis auf \glq wenige\grq\ Ausnahmen. Ob es für alle $\Omega$ gilt, ist
unbekannt. Referenzen für diese Resultate sind \cite{Cha:FPVR} -- für Mannigfaltigkeiten ohne Rand -- und \cite{GuiMel:PSFMWB} im technisch schwierigeren Fall mit Rand, wobei hier immer angenommen wird, dass es keine Ecken gibt. Ecken bereiten noch größere Probleme. Ob umgekehrt das Längenspektrum die Eigenwerte bestimmt, ist unbekannt.}

\bigskip
Was haben Obertöne mit Billiardkugelbahnen zu tun? Vielleicht erinnern Sie sich an die Welle-Teilchen-Dualität im Physik-Unterricht der Schule. Licht (oder Elektronen oder...) ist als Welle zu beschreiben -- ähnlich wie Schall, daher die Beziehung zu Obertönen -- und gleichzeitig mittels Partikeln, die sich geradlinig fortbewegen -- wie eine Billiardkugel!

Mathematisch werden Wellen durch eine Differentialgleichung beschrieben, die \emph{Wellengleichung}:\\[2ex]
\setlength{\fboxrule}{1pt} \setlength{\fboxsep}{2mm}
\framebox{\parbox{12cm}{\begin{align}  \label{eq-wgl}
\left( \frac{\partial}{\partial t^2} - \Delta\right) w(t,x) & = 0, \qquad
t\in \R, \; x\in \Omega \notag\\
w(0,x) & =f(x)\notag\\[-1ex]
 & \tag{WGl}\\[-1ex]
\frac{\partial}{\partial t}\; w(0,x) & = 0\notag\\[1ex]
w(t,x) & =0, \qquad t\in \R, \; x\in \partial \Omega \ (=\text{Rand von }\Omega).\notag
\end{align}}
}
\bigskip

Das sieht sehr kompliziert aus, hat aber eine einfache konkrete Bedeutung:
Betrachtet man $\Omega$ als Trommelfell und bringt man dies zum Zeitpunkt
$t=0$ in die Form $f$ (z.~B. erhält man ein ``hutförmiges'' $f$, etwa $f(x)
= 1-|x|$, wenn man bei einer kreisförmigen Trommel -- $\Omega = \{x: \; |x|
< 1\}$ -- die Mitte nach oben drückt) und lässt dann los, so beschreibt die
Lösung $w(t,\cdot)$ die Form der Trommel zum Zeitpunkt $t$.  Die erste Gleichung beschreibt die schwingungserzeugenden Kräfte (idealisiert, d.~h. ohne Reibung und in linearer Näherung), die zweite und dritte die Anfangsposition und die vierte die Tatsache, dass die Trommel am Rand fest eingespannt ist.

Man kann zeigen, dass \eqref{eq-wgl} für jedes stetige $f$ eine
eindeutige Lösung $w$ besitzt. Für jedes $t$ schreibe $W_t(f) =
w(t,\cdot)$ für den \glq Lösungsoperator zur Zeit $t$\grq.

Wir wenden nun die Spurformel auf $W_t$ an (für beliebiges festes \(t\)). Die Eigenwerte von $W_t$ lassen
sich leicht aus denen von $\Delta$ gewinnen. Gilt nämlich $-\Delta u(x) =
\mu u(x)$, $u_{|\partial \Omega} = 0$, so ist
\[
   w(t,x) = \cos (\sqrt{\mu} t) \cdot u(x)
\]
eine Lösung von \eqref{eq-wgl} mit \(f=u\), wie man sofort nachrechnet. Also \(W_t(u) = \cos (\sqrt \mu t) u\). Sind also $0 <
\mu_1 \leq \mu_2 \leq \ldots$ die Eigenwerte von $-\Delta$, so sind $\cos
\sqrt{\mu_1} t, \cos \sqrt{\mu_2} t,\ldots, $ die Eigenwerte von $W_t$. Die
linke Seite der Spurformel für $W_t$ ist also
\begin{equation} \label{eq24}
   \sum^\infty_{k=1} \cos \sqrt{\mu_k} t\, .
\end{equation}
Dass dies für die meisten $t$ (z.~B. für $t=0$) divergiert, soll uns
zunächst nicht kümmern.

Um die rechte Seite der Spurformel auszuwerten, brauchen wir eine alternative
Methode, eine Funktion $G_t(x,y)$ zu finden mit
\[
   (W_tf)(x) = \int_{\Omega} G_t(x,y)\, f(y)\, dy\, , \qquad x \in \Omega
\]
für alle $f$ (dies ist der zweite Schritt in der Strategie in Abschnitt \ref{secneublick}.

In diesem Artikel können wir hierzu nur die Idee skizzieren, die man am besten
anhand der physikalischen Vorstellung der Wellenausbreitung versteht: Wir
stellen uns nun $\Omega$ als Teich vor. Zum Zeitpunkt $t=0$ werfen wir am Ort
$y$ einen Stein ins Wasser. Danach breitet sich eine Wellenfront
kreisförmig um $y$ aus. $G_t(x,y)$ ist die Höhe der Welle zum Zeitpunkt $t$ am Ort
$x$.\footnote{Denn das Hineinwerfen des Steins beim Punkt \(y\) kann, idealisiert, durch \(f=\delta_y\), die bei \(y\) konzentrierte Delta-Funktion -- genauer $\delta-$Distribution, siehe Fußnote \ref{fndistr} --, beschrieben werden, und es ist \((W_t (\delta_y))(x) = G_t(x,y)\). Genau genommen
entspricht dieses Experiment der Version der Wellengleichung, wo $u(0,x)
= 0$, $u_t(0,x) = f$ ist, aber das macht für die Methode keinen
Unterschied.}
Wie kann man die Wellenfront geometrisch beschreiben? Die Wellenfront besteht aus den Punkten, die man von \(y\) aus entlang einer geraden Linie der Länge \(t\) erreichen kann. Sobald diese Linie den Rand erreicht, wird sie nach dem Reflektionsgesetz reflektiert. Damit ergeben sich für kleine $t$ kreisförmige Wellenfronten (vom Radius $t$), und sobald der Kreis den Rand erreicht, ergeben sich kompliziertere Muster.

Also:
\begin{center}
\begin{tabular}{lcp{6cm}}
Wellenfront nach der Zeit $t$ &  =  & Endpunkte der Billiardkugelbahnen der
Länge $t$, die in $y$ beginnen.
\end{tabular}
\end{center}
Damit folgt: $G_t(x,x)$ ist genau dann besonders
groß, wenn es eine Billiardkugelbahn der Länge $t$ gibt, die in $x$ startet
und in $x$ endet.

Kommen wir nun zur Spurformel für $W_t$. Mit \eqref{eq24} lautet sie
\begin{equation} \label{eq25}
   \sum^{\infty}_{k=1} \cos  \sqrt{\mu_k} t\; = \int_\Omega G_t (x,x)\, dx
    \, , \qquad t\in \R\, .
\end{equation}
Nach diesen Überlegungen wird die rechte Seite nur für diejenigen Werte $t$
besonders groß, für die es eine Billiardkugelbahn der Länge $t$ von $x$ nach
$x$ gibt, für mindestens ein $x$. Eine genauere Analyse ergibt,
dass nach der Integration über $x$ nur von solchen Bahnen ein wesentlicher Beitrag
übrigbleibt, die auch in der Startrichtung wieder zu $x$ zurückkehren, die
also geschlossen sind.

Wie fasst man all dies mathematisch, was heißt ``besonders groß'', wie
steht es mit der Divergenz der linken Seite von \eqref{eq25}?

Eine elegante Antwort bietet die Distributionentheorie: Wir können hier
nicht die Definition von Distributionen geben. Man kann sie sich grob als
Funktionen vorstellen, die an einigen Stellen in einem sehr präzisen Sinn
unendlich sind.\footnote{\label{fndistr}Und zwar derart, dass auch isolierte
$\infty$-Punkte bei Integration einen Beitrag liefern. Standardbeispiel:
Die $\delta$-Distribution (auf $\R$), die für $x\not= 0$ gleich Null ist
und bei $x=0$ in einer Weise $\infty$ ist, dass ``$\int\delta (x)\,
\varphi(x)\, dx = \varphi(0)$'' für beliebige glatte Funktionen $\varphi$
ist.}
Man nennt sie an diesen Stellen \emph{singulär}. Nun ist $G_t(x,y)$ in
Wirklichkeit eine Distribution, und die Reihe $\sum^\infty_{k=1} \cos
\sqrt{\mu_kt}$ konvergiert im Sinne der Distributionen, d.~h. die Summe ist
eine Distribution. Ersetzt man noch
``besonders groß'' durch ``singulär'', so werden die Überlegungen oben zu
einer mathematischen (und korrekten) Aussage: Die Menge der $t>0$, an denen
die rechte Seite der Spurformel \eqref{eq25} singulär ist, ist gleich der
Menge der Längen geschlossener Billiardkugelbahnen.

Da diese Menge gleichzeitig durch die linke Seite der Spurformel bestimmt
ist, also durch die $\mu_k$, ergibt sich die Behauptung am Anfang dieses
Kapitels.

\section{Spurformel und Theta-Funktion}
Wir betrachten hier eine dritte frappierende Anwendung der Spurformel. Wir gehen ähnlich vor wie im vorhergehenden Abschnitt, aber statt der Wellengleichung betrachten wir die
\emph{Wärmeleitungsgleichung}:
\begin{align} \label{eq-wlo}
   \left(\frac{\partial}{\partial t} - \frac{\partial^2}{\partial x^2}
   \right) u(t,x) & = 0 \qquad t > 0 \notag\\
   u(0,x) & = f(x) \tag{WLGl}\\
   u(t,x+1)&=u(t,x). \notag
\end{align}
Hierbei ist immer $t>0$ und $x\in \R$. Die letzte Bedingung ist als Variante zur Randbedingung $u(0)=u(1)=0$ in \eqref{eq8} zu verstehen. Sie bedeutet, dass $u$ bzgl. $x$ $1$-periodisch ist; natürlich muss $f$ auch $1$-periodisch sein. Auf die physikalische Bedeutung von (WLGl) gehen wir hier nicht ein (siehe z.B.\ \cite{CouHil:MMP}). Die Spurformel liefert hier die höchst bemerkenswerte Identität (für \(t>0\) beliebig)
\begin{equation} \label{eq-wl}
   \sum^\infty_{k=-\infty} e^{-t\cdot 4\pi^2k^2} = \frac{1}{\sqrt{4\pi t}}
     \; \sum^\infty_{\ell=-\infty} e^{-\frac{1}{t} \cdot \frac{\ell^2}{4}}
\end{equation}
Dies wird üblicherweise $\theta(s) = \frac{1}{\sqrt{s}} \; \theta
\left(\frac{1}{s}\right)$ mit $\theta(s) = \sum^\infty_{k=-\infty}
e^{-s\pi k^2}$ geschrieben (ersetze $s = 4 t\pi$). Diese
\emph{Transformationsformel der Theta-Funktion} ist in verschiedenen
Bereichen der Mathematik fundamental, z.~B. bildet sie die Basis der
Theorie der Modulformen, die unter anderem in der Zahlentheorie spannende Anwendungen hat.

Man erhält \eqref{eq-wl} wie folgt: Wiederum kann man zeigen, dass (WLGl) für jedes stetige $1$-periodische $f$ eine eindeutige Lösung hat.\footnote{Die Existenz folgt aus der weiter unten angegebenen Formel. Die Eindeutigkeit zeigt man so: Ist $u$ eine Lösung von (WLGl), so gilt $\frac d{dt} \int_0^1 u^2(t,x)\,dx = 2\int_0^1u(t,x)\partial_t u(t,x)\,dx = 2\int_0^1 u(t,x)\partial_x^2u(t,x)\,dx = -2\int_0^1 (\partial_x u(t,x))^2\,dx\leq 0$, wobei im letzten Schritt partiell integriert wurde unter Verwendung der periodischen Randbedingung. Sind nun $u_1,u_2$ zwei Lösungen mit derselben Anfangsbedingung, so ist auch $u=u_1-u_2$ Lösung mit Anfangsbedingung $0$. Da $\int_0^1 u^2(t,x)\,dx$ als Funktion von $t$ monoton fällt, immer nicht-negativ ist und bei $t=0$ gleich Null ist, muss es konstant gleich Null sein, also $u\equiv 0$, also $u_1\equiv u_2$.} Zu $t>0$ sei $P_t: f \longmapsto
u(t,\cdot)$ der Operator, der die Anfangsfunktion $f$ in die Lösung zur
Zeit $t$ überführt. Die Eigenwerte von $P_t$ lassen sich leicht aus denen von
$-\frac{d^2}{dx^2}$ gewinnen: Ist $-\frac{d^2}{dx^2}\; f = f$, so ist
offenbar
\[
   u(t,x) = e^{-t\mu}  f(x)
\]
eine Lösung von \eqref{eq-wlo}, also ist \(e^{-t\mu}\) ein Eigenwert von \(P_t\). Ähnlich wie in Abschnitt \ref{secewberechnung} zeigt man, dass die $1$-periodischen Eigenfunktionen von
$-\frac{d^2}{dx^2}$ genau $f_k(x) = e^{2\pi ikx}$ sind, mit Eigenwerten  $\mu_k =
4\pi^2k^2$. Also ist die Spur von \(P_t\) gerade die linke Seite von \eqref{eq-wl}.
Konvergenzprobleme gibt es hier nicht.

Jetzt brauchen wir noch eine alternative Lösungsmethode für \eqref{eq-wlo}.
Betrachtet man die Wärmeleitungsgleichung ohne die Periodizitätsbedingung an $u$ und $f$, so kann man eine Lösung explizit hinschreiben:
\[
   u(t,x) = \int^\infty_{-\infty} K_t (x,y) \, f(y)\, dy
\]
mit
\[
   K_{t}(x,y) = \frac{1}{\sqrt{4\pi t}}  \; e^{-{|x-y|^2}/{4 t}}
   \, ,
\]
wie man durch Einsetzen (unter Verwendung von $\int^\infty _{-\infty}
e^{-x^2}dx = \sqrt{\pi}$) direkt nachprüft.%
\footnote{Diese Formel kann man z.~B. mittels der Fouriertransformation, aber auch
durch Zurückführung auf eine gewöhnliche Differentialgleichung mittels
Homogenitätsüberlegungen, herleiten.}

Ist nun $f$ $1$-periodisch, so sieht man leicht, dass dann
\[
   (P_t f)(x) = u(t,x) = \int^1_0  K^{\text{per}}_t(x,y) f(y)\, dy
\]
mit
\[
   K^{\text{per}}_t(x,y) = \sum^\infty_{\ell=-\infty} K_t(x,y+\ell)
\]
eine $1$-periodische Lösung von (WLGl) ist. Indem man $x=y$ setzt und über $x$ von 0 bis 1 integriert,
erhält man die rechte Seite von \eqref{eq-wl}.

\section{Schlussbemerkungen}
Wir haben an drei Beispielen illustriert, welche Welten sich eröffnen, wenn man die einfache Spurformel für Matrizen verallgemeinert. In jedem Fall stand eine Differentialgleichung im Hintergrund, und wesentlich war, dass man neben der spektralen eine direkte Lösungsmethode für diese Gleichung hatte. Der Laplace-Operator $P=-\Delta$ (mit Dirichlet-Randbedingungen im ersten und zweiten Fall) trat in allen Beispielen auf.  Bei genauerem Hinsehen erkennt man, dass im ersten Beispiel ($-u''=f$) die Spurformel auf den Operator \(P^{-1}\), im zweiten (Wellengleichung)
auf den Operator \(\cos t\sqrt P\) und im dritten Beispiel auf den Operator \(e^{-tP}\) angewendet wurde. Führt man ein ähnliches Programm mit dem Schrödinger-Operator \(e^{-i\frac t\hbar P}\) durch, erhält man die Gutzwiller-Spurformel, die die Korrespondenz zwischen klassischer und Quantenmechanik zum Ausdruck bringt (\cite{Gut:POCQC}). Eine weitere berühmte Spurformel ist die von Selberg im Kontext der Zahlentheorie (\cite{Hej:STFP}). Diese erhält man durch Betrachtung des Wärmeleitungsoperators auf Flächen negativer Krümmung (statt des Gebietes \(\Omega\)).

Natürlich gibt es für \(\sum_k \frac1 {k^2} =  \frac{\pi^2}6\)  auch andere Beweise, ganz ohne Differentialgleichungen (siehe z.B.\ \cite{Cha:EZ} für eine hübsche Sammlung), ebenso für die Transformationsformel der Theta-Funktion \eqref{eq-wl}.
Die Spurformel gibt jedoch ein faszinierendes Werkzeug an die Hand, das so
diverse Dinge wie Welle-Teilchen-Dualität, Zahlentheorie und
$\sum \frac{1}{k^2} = \frac{\pi^2}{6}$ in Beziehung setzt. Vielleicht
finden Sie ja noch weitere  interessante Anwendungen?

\theendnotes


\def\cprime{$'$}

\end{document}